\documentclass[11pt,number,sort&compress]{elsarticle}

\usepackage{amsmath}
\usepackage{mathtools}
\usepackage{mathrsfs}
\usepackage{empheq}
\usepackage{amsfonts}
\usepackage{setspace}
\usepackage{cleveref}
\usepackage{empheq}
\usepackage{appendix}
\usepackage{color}
\usepackage{lineno}

\def\L{\left}
\def\R{\right}

\newcommand{\Dc}[2]{\L(\prescript{C}{}{\mathbf{D}}_{0^+}^{#2} #1 \R)\L(t\R)} 
\newcommand{\IRL}[2]{\L({\mathbf{I}}_{0^+}^{#2} #1 \R)\L(t\R)} 

\begin{document}

\begin{frontmatter}

\title{Fractional differential equations solved by using Mellin transform}

\author[HW]{Salvatore Butera\corref{cor1}}
\ead{sb469@hw.ac.uk}
\cortext[cor1]{Corresponding author}
\address[HW]{SUPA, Institute of Photonics and Quantum Sciences, Heriot-Watt University, Edinburgh, EH14 4AS, United Kingdom.}

\author[unipa]{Mario Di Paola}
\ead{mario.dipaola@unipa.it}
\address[unipa]{Dipartimento di Ingegneria Civile Ambientale, Aerospaziale e dei Materiali (DICAM), Universit\`{a} degli Studi di Palermo, Viale delle Scienze,Ed.8, 90128 - Palermo, Italy}


\begin{abstract}
In this paper, the solution of the multi-order differential equations, by using Mellin Transform, is proposed. It is shown that the problem related to the shift of the real part of the argument of the transformed function, arising when the Mellin integral operates on the fractional derivatives, may be overcame. Then, the solution may be found for any fractional differential equation involving multi-order fractional derivatives (or integrals). The solution is found in the Mellin domain, by solving a linear set of algebraic equations, whose inverse transform gives the solution of the fractional differential equation at hands.
\end{abstract}


\end{frontmatter}

\section{Introduction}
It is widely recognized that the memory and hereditary properties of various materials and processes in  electrical circuits, biology, biomechanics, etc. such as viscoelasticity, is well predicted by using fractional differential operators. Such operators are the generalization, to real (or complex) order, of the classical derivatives and integrals (see e.g. \cite{miller1993,mainardi1997,oldham2006}). 
Conversely to the locally defined classical derivatives, the powerful of the fractional operators in describing the time evolution of many physical processes and, in general, in modelling the dynamics of complex systems, is due to the long memory characteristics inherent to these operators. Indeed, dealing with systems characterized by a power type non local interaction, or by a non-Markovian power law time memory, into which the complexity of the dynamics usually manifest itself, fractional differential equations naturally arise in the relative mathematical models \citep{Tarasov2008,Cottone2009,DiPaola2008}. As a proof of their powerful in describing nature, theoretical research on these operators experienced an exceptional boost in the last few decades, and applications can now be found in various fields of natural sciences. Examples are in electrical circuits \cite{Ala2013}, in anomalous transport and diffusion processes in complex media \cite{Schneider1989, Isichenko1992, Zaslavsky2002, Metzler2001}, in material sciences \cite{MainardiSurvey,Koeller1984,Bagley1986}, in biology \citep{Xu2009,Arafa2012,Ahmed2012} and biomechanics \citep{Deseri2013,Magin2010a,Magin2010b}, and in many other branches of physics and engineering \cite{podlubny1998,kilbas2006,mainardi2010}.

Various methods for the solution of differential equation of fractional order are available in literature, including Laplace method \cite{podlubny1998,Duan2004}, Gr\"{u}nwald-Letnikov method \cite{podlubny1998, Duan2008}, Adomian method \cite{Wazwaz1996} and several others \cite{podlubny1998,He1998,Wang2007,Song2007,Xu2009,Saadatmandi2010}. Some attempts to use Mellin Transform and related concepts, have been presented (see e.g. \cite{podlubny1998}) in order to solve particular classes of fractional differential equations.

In this paper, a general method of solution for Initial Value Problems (IVP), involving fractional derivatives, is presented by using the Mellin Transform of complex order $\gamma=\rho+i \eta$. The method takes advantage of the fact that the discretized version of the inverse Mellin transform may be seen, in logarithmic temporal scale, as a Fourier series and, in time domain, as a complex Taylor series with coefficients depending on the fractional integrals in zero. The restitution of the function is independent on the value of $\rho$ used to evaluate the discretized inverse Mellin transform, provided it belongs to its so called \textit{fundamental strip}. With these informations in mind, in this paper, a method that allow us to relate the value of the Mellin Transform for different values of $\rho$ belonging to the fundamental strip, is presented. This is the main key to solve multi-order fractional differential equations in a very easy and direct way. The method is versatile and easy to implement in computer programs.

The paper is organized as follows: in the next section, the need to handle with fractional differential equation is presented, with a relevant example in which fractional differential equations appear and whose classical solution is written in terms of Mittag-Leffler series expansion; in section 3, the Mellin Transform and related concepts are highlighted; in section 4, the solution of the fractional differential equation is presented, along with some applications in sect. 5. In appendix some few basic elements on fractional calculus are reported for completeness sake's.

\section{Fractional differential equations hereditarieness}
In this section, the relevant example of viscoelastic materials is presented in order to show the importance of the fractional calculus for many physical and engineering problems.

The linear viscoelastic problem is ruled by two different but interconnected functions: i) the creep function, labelled as $J(t)$, that is the strain history for an imposed stress history $\sigma(t)=U(t)$ (unit step); ii) the relaxation function, labelled ad $G(t)$, that is the stress history for an imposed strain history $\epsilon(t)=U(t)$. In linear viscoelasticity, the Boltzmann superposition principle holds, so that
\vspace{0mm}
\begin{subequations}
\label{eq:1}
\begin{equation}
    \sigma(t)=\int_0^t{G(t-\tau)\dot{\epsilon}(\tau)d\tau}
\label{eq:1a}
\end{equation}
\begin{equation}
       \epsilon(t)=\int_0^t{J(t-\tau)\dot{\sigma}(\tau)d\tau}
\label{eq:1b}
\end{equation}
\end{subequations}
\vspace{0mm}
Eqs.\eqref{eq:1}, valid for quiescent systems in $t\leq 0$, suggest that relaxation and creep functions play the role of kernels, in \eqref{eq:1a} and \eqref{eq:1b} respectively. By using the Laplace transform of eqs.\eqref{eq:1}, the following fundamental relationship 
\vspace{0mm}
\begin{equation}
    \hat{J}(s)\hat{G}(s)=1/s^2
\label{eq:2}
\end{equation}
\vspace{0mm}
\noindent is easily derived, where $\hat{J}(s)$ and $\hat{G}(s)$ are the Laplace transform of $J(t)$ and $G(t)$ respectively.\\
On the other hand, Nutting \cite{Nutting1921}, by means of experimental tests performed on various materials like rubber, ceramics, etc., showed that in general the relaxation function may be written in the form
\vspace{0mm}
\begin{equation}
    G(t)=\frac{c_\alpha}{\Gamma\L(1-\gamma\R)} \, t^{-\alpha} \qquad (0\leq\alpha\leq 1)
\label{eq:3}
\end{equation}
\vspace{0mm}
and, in virtue of eq.\eqref{eq:2}, the corresponding creep function is 
\vspace{0mm}
\begin{equation}
    J(t)=\frac{1}{c_\alpha\Gamma\L(1+\gamma\R)} \, t^{\alpha} \qquad (0\leq\alpha\leq 1)
\label{eq:4}
\end{equation}
\vspace{0mm}
\noindent where $\Gamma\L(\cdot\R)$ is the Euler Gamma function, $c_\alpha$ and $\alpha$ are characteristic coefficients of the material at hands.\\
As we insert eqs. \eqref{eq:3} and \eqref{eq:4} in eqs.\eqref{eq:1a} and \eqref{eq:1b}, respectively we get
\vspace{0mm}
\begin{subequations}
\label{eq:5}
    \begin{align}
      \sigma(t)&=c_\alpha \Dc{\epsilon}{\alpha} \label{eq:5a} \\
       \epsilon(t)&=\frac{1}{c_\alpha}\IRL{\sigma}{\alpha} \label{eq:5b}
    \end{align}
  \end{subequations}
\vspace{0mm}
\noindent where $\Dc{\epsilon}{\alpha}$ and $\IRL{\sigma}{\alpha}$ are the Caputo's functional derivative \cite{Gerasimov1948} and the Riemann-Liouville fractional integral, respectively (see appendix). From eqs.\eqref{eq:5}, some considerations may be drawn: i) the viscoelastic  constitutive law is ruled, in its direct and inverse form, by fractional operators (derivative and integral) of the same order $(\alpha)$. ii) if $\alpha=0$, then the elastic constitutive law is recovered while, if $\alpha=1$, the Newton-Petrov constitutive law of pure fluid appears. It follows that the constitutive law of a viscoelastic material has an intermediate behaviour between pure fluid and pure elastic solid. iii) for quiescent systems for $t\leq 0$, the Caputo's fractional derivative coalesces with the Riemann-Liouville fractional derivative, and the two operators in eqs. \eqref{eq:5a} and \eqref{eq:5b} are the inverse each another.

In order to capture different behaviours for more complex systems like bones \cite{Deseri2013}, bitumen \cite{Grzesikiewicz2013}, and so on, the constitutive law has to be modified by inserting others fractional order operators, to obtain differential equations of the kind
\vspace{0mm}
\begin{equation}
    \sigma(t)=\sum_{k=0}^n{c_k^\alpha \Dc{\epsilon}{\alpha_k}}
\label{eq:6}
\end{equation}
\vspace{0mm}
\noindent With this representation of the constitutive law, the solution in term of the strain for an imposed stress history, shows some troubles. In order to show this, let us start with the simple case
\vspace{0mm}
\begin{subequations}
\label{eq:7}
\begin{empheq}[left=\empheqlbrace]{align}
     \sigma(t)&={c_0\epsilon(t)+c_1^{\alpha}\Dc{\epsilon}{\alpha}} \label{eq:7a}\\
     \epsilon(0)&=0 \label{eq:7b}
\end{empheq}
\end{subequations}
\vspace{0mm}
\noindent that is the Kelvin-Voigt element in which the dashpot is substituted by a fractional viscoelastic element. By making the Laplace transform of eq.\eqref{eq:7}, we obtain

\begin{equation}
	\hat{\epsilon}(s)\L[c_0+c_1^{(\alpha)}s^\alpha\R]-s^{\alpha-1}\epsilon(0)=\hat{\sigma}(s)
\label{eq:8}
\end{equation}
\vspace{0mm}
\noindent and, by inserting the initial condition \eqref{eq:7b} in eq.\eqref{eq:8}, we get
\vspace{0mm}
\begin{equation}
	\hat{\epsilon}(s)=\frac{\hat{\sigma}(s)}{c_0+c_1^{\alpha} s^\alpha}
\label{eq:9}
\end{equation}
\vspace{0mm}
\noindent The inverse Laplace transform is then \cite{podlubny1998}
\vspace{0mm}
\begin{equation}
\begin{split}
    \epsilon(t)&=\frac{1}{c_1^\alpha}\int_0^t{\L(t-\tau\R)^{\alpha-1}E_{\alpha,\alpha}\L(-\frac{c_0}{c_1^\alpha} \L(t-\tau\R)^\alpha\R)\sigma(\tau)\, d\tau}\\
    &=\frac{1}{c_1^\alpha}\sum_{k=0}^\infty{\L[\L(-\frac{c_0}{c_1^\alpha}\R)^{k}\IRL{\sigma}{\alpha(k+1)}\R]}
    \end{split}
\label{eq:10}
\end{equation}
\vspace{0mm}
\noindent where the two parameter Mittag-Leffler function \cite{podlubny1998} is defined as
\vspace{0mm}
\begin{equation}
	E_{\alpha,\beta}(t)=\sum_{k=0}^{\infty}{\frac{t^k}{\Gamma\L(\alpha+k\beta\R)}} \qquad \L(\alpha>0,\; \beta>0\R)
\label{eq:11}
\end{equation}
\vspace{0mm}
\noindent In spite eq.\eqref{eq:10} is correct from a mathematical point of view, it presents a problem connected with its definition in terms of Mittag-Leffler function, that is the truncation of the series pathologically diverges at infinity.\\
The method proposed in the next section will overcame this problem by using the Mellin transform.

\section{Mellin Transform}
The Mellin transform of a given function $x(t)$, defined in the range $0\leq t\leq\infty$, is given as
\vspace{0mm}
\begin{equation}
    X\L(\gamma\R)=\mathcal{M}\L\{x(t);\gamma\R\}\equiv\int_0^\infty{t^{\gamma-1}\,x(t)\,dt}; \qquad \gamma=\rho+i \eta
\label{eq:12}
\end{equation}
\vspace{0mm}
\noindent and its inverse transform is
\vspace{0mm}
\begin{equation}
    x\L(t\R)=\mathcal{M}^{-1}\L\{X(\gamma);t\R\}\equiv \frac{1}{2\pi}\int_{-\infty}^{\infty}{X(\gamma)\,t^{-\gamma}\, d\eta};\qquad (t>0)
\label{eq:13}
\end{equation}
\vspace{0mm}
\noindent Condition of existence of \eqref{eq:12} and \eqref{eq:13} is $-p<\rho<-q$, where $p$ and $q$ refer to the asymptotic behaviour of the function at $t\to 0$ and $t\to\infty$, respectively:
\vspace{0mm}
\begin{equation}
	x(t)\sim t^p \quad (t\to 0)\; ;\quad x(t)\sim t^q \quad (t\to \infty)
\label{eq:14}
\end{equation}
\vspace{0mm}
\noindent Equation \eqref{eq:13} may be discretized in the form
\vspace{0mm}
\begin{equation}
    x(t)\simeq \frac{\Delta\eta}{2\pi}\sum_{k=-m}^{m}{X(\gamma_k) t^{-\gamma_k}}; \qquad \gamma_k=\rho+i k \Delta \eta
\label{eq:15}
\end{equation}
\vspace{0mm}
\noindent where $m\Delta\eta=\bar{\eta}$ is a cutoff value along to the immaginary axis, selected in such a way that the contribution of the $m+n-\text{th}$ terms  in the summation \eqref{eq:15} do not produce sensible variation on $x(t)$. It is to be noticed that the Mellin transform is holomorphic into the fundamental strip, and $X\L(\gamma_k\R)=X^*\L(\gamma_{-k}\R)$. Moreover, the result of the summation is independent of the value $\rho$, provided it is selected into the fundamental strip. It is worth to stress that eq.\eqref{eq:15} may be also rewritten in the form
\vspace{0mm}
\begin{equation}
    x(t)\simeq t^{-\rho}\L\{\frac{A_0}{2b}+\frac{1}{b}\sum_{k=1}^{m}{\L[A_k \cos\L(\frac{k\pi}{b}\text{ln}\, t\R)+B_k \sin\L(\frac{k\pi}{b}\text{ln}\, t\R)\R]}\R\}
\label{eq:16}
\end{equation}
\vspace{0mm}
\noindent with $b=\pi/\Delta\eta$ and $A_k=\text{Re}\L[X{\L(\gamma_k\R)}\R]$ and $B_k=\text{Im}\L[X{\L(\gamma_k\R)}\R]$. From eq.\eqref{eq:16}, we recognize that the term in parenthesis is a Fourier series in logarithmic temporal scale. It is easy to show that the series has a the time dependent period: $T=\big (\exp(2b)-1\big )\,t$. This feature, along with the factor $t^{-\rho}$ in front of the series, reveals a hidden self-similar structure of the whole function $x(t)$, with a fractal Hurst dimension $H=\rho$. This means that exist $a>0$, such that $x(at)=a^H\;x(t)$ \cite{gennady1994}. We demonstrate now the validity of this sentence for the particular value $a=\exp(2b)$. Denoting with $\tilde{x}(t)$ the series in eq.\eqref{eq:16}, and taking into account the expression of the period $T$ given above, we have that $\tilde{x}(t)=\tilde{x}(t+T)=\tilde{x}(\exp(2b)t)=\tilde{x}(at)$. From this result follows the relation $x(at)=a^{-\rho}t^{-\rho}\tilde{x}(at)=a^{-\rho}t^{-\rho}\tilde{x}(t)=a^{-\rho}x(t)$, that demonstrates the statement given above.

Let us now suppose that we are involved in finding the solution of the following fractional differential equation of order $n-1<\alpha<n$, with $n \in \mathbb{N}$
\vspace{0mm}
\begin{subequations}
\label{eq:17}
\begin{empheq}[left=\empheqlbrace]{align}
     &\Dc{x}{\alpha}+\lambda x\L(t\R)=f(t) \label{eq:17a}\\
     &\;\;x(0)=0\; ; \quad ... \quad x^{(n-1)}(0)=0 \label{eq:17b}
\end{empheq}
\end{subequations}
\vspace{0mm}
\noindent that generalize, to arbitrary order, the Cauchy problem defined in eqs.\eqref{eq:7} for $0<\alpha<1$. The Mellin transform of such fractional differential equation is given as (see, \cite{kilbas2006})
\vspace{0mm}
\begin{equation}
\begin{split}
    \mathcal{M} \L\{\Dc{x}{\alpha};\gamma\R\}=&\sum_{k=0}^{n-1}{\frac{\Gamma\L(1-\gamma-n+\alpha+k\R)}{\Gamma\L(1-\gamma\R)}\L[x^{(n-k-1)}(t)\, t^{\gamma+n-\alpha-k-1}\R]_0^\infty}\\
    &+\frac{\Gamma\L(1-\gamma+\alpha\R)}{\Gamma(1-\gamma)}X\L(\gamma-\alpha\R)
\label{eq:18}
\end{split}
\end{equation}
\vspace{0mm}
\noindent In the case $\lambda>0$, the system is stable so that, taking into account the initial conditions \eqref{eq:17b}, and assuming $f(t)$ zero from some value $t=t_{\text{max}}$, then $x(t)\to0$ for $t\to\infty$, along with all its derivatives. The assumption $\rho<\{\alpha\}$, with $\{\alpha\}$ the fractional part of $\alpha$, guarantees that the first term in the left-hand side of eq.\eqref{eq:18} vanishes, so the Mellin transform of eq.\eqref{eq:17a} takes the form
\vspace{0mm}
\begin{equation}
    C\L(\gamma,\alpha\R)X\L(\gamma-\alpha\R)+\lambda X(\gamma)=F(\gamma)
\label{eq:19}
\end{equation}
\vspace{0mm}
\noindent being $F(\gamma)=\mathcal{M}\L\{f(t);\gamma\R\}$ and
\vspace{0mm}
\begin{equation}
	C\L(\gamma,\alpha\R)=\frac{\Gamma\L(1-\gamma+\alpha\R)}{\Gamma(1-\gamma)}
\label{eq:20}
\end{equation}
\vspace{0mm}
\noindent From eq.\eqref{eq:19}, we recognise that the solution in terms of Mellin transform may not be pursued since $X(\gamma-\alpha)$ is not related at $X(\gamma)$. For this reason, only particular classes of fractional differential equations, in which $\lambda=\lambda(t)$ and properly selected, may be solved (see \cite{podlubny1998}).
\vspace{0mm}
\noindent Recently, Di Paola showed that $X(\gamma-\alpha)$ and $X(\gamma)$ may be related each another in a very simple way. This relation takes full advantage from the fact that, according to eq.\eqref{eq:15}, $x(t)$ is indipendent to the value of $\rho$ selected into the fundamental strip. Then we may write
\vspace{0mm}
\begin{equation}
    x(t)\simeq \frac{\Delta\eta}{2\pi}\sum_{k=-m}^{m}{X(\gamma_k-\alpha) t^{-(\gamma_k-\alpha)}}= \frac{\Delta\eta}{2\pi}\sum_{k=-m}^{m}{X(\gamma_k) t^{-\gamma_k}}
\label{eq:21}
\end{equation}
\vspace{0mm}
\noindent that remain valid if $\rho-\alpha$ (and of course $\rho$) belong to the fundamental strip. Let us suppose that the $X\L(\gamma_k\R)$ is known and let us evaluate $X\L(\gamma_k-\alpha\R)$. To this aim, we rewrite eq.\eqref{eq:21} in the form
\vspace{0mm}
\begin{equation}
    t^{-\frac{1}{2}}\sum_{k=-m}^{m}{X(\gamma_k-\alpha) \exp\L({-i \frac{k\pi}{b}\ln t}\R)}\simeq t^{-\L(\alpha+\frac{1}{2}\R)}\sum_{k=-m}^{m}{X(\gamma_k) \exp\L({-i \frac{k\pi}{b}\ln t}\R)}
\label{eq:22}
\end{equation}
\vspace{0mm}
\noindent and we require that the modulus of the difference, integrated over a proper domain, is minimum, that is
\vspace{0mm}
\begin{multline}
    \int_{t_1}^{t_2}\frac{1}{t}\L[\sum_{k=-m}^{m}{X(\gamma_k-\alpha) \exp\L({-i \frac{k\pi}{b}\ln t}\R)}-t^{-\alpha}\sum_{k=-m}^{m}{X(\gamma_k) \exp\L({-i \frac{k\pi}{b}\ln t}\R)}\R]\\
    \times\L[\vphantom{\sum_{k=-m}^{m}} \text{c.c.}\R]dt=\underset{X(\gamma_k-\alpha)}{\text{min}}
\label{eq:23}
\end{multline}
\vspace{0mm}
\noindent where $\L[\text{c.c.}\R]$ stands for complex conjugate.\\
In principle, $t_1$ and $t_2$ should be $t_1=0$ and $t_2=\infty$ but, because in $t=0$ singularities appears, we select $t_1=\exp\L({-b}\R)$ and $t_2=\exp\L({b}\R)$. In this way we have that the overlapping of the response, evaluated for two values of the fundamental strip $\L(\rho\,\text{and}\, \rho-\alpha\R)$, is guaranteed in a very large interval. Now we make the change of variable
\vspace{0mm}
\begin{equation}
    \ln t=\xi\,; \quad \frac{dt}{t}=d\xi\,; \quad\ln t_1=-b\,; \quad \ln t_2=b
\label{eq:24}
\end{equation}
\vspace{0mm}
\noindent By performing variations respect to the $X(\gamma_k-\alpha)$ $(k=-m,...,m)$ and taking advantage from the orthogonality of $\exp\L({-i\frac{k\pi}{b}\xi}\R)$ on $[-b, b]$, we get
\vspace{0mm}
\begin{equation}
    X\L(\gamma_s-\alpha\R)=\frac{1}{2b}\sum_{k=-m}^{m}{X(\gamma_k)a_{sk}(\alpha)}
\label{eq:25}
\end{equation}
\vspace{0mm}
\noindent where
\vspace{0mm}
\begin{equation}
    a_{sk}(\alpha)=\int_{-b}^b{\exp\L({-\L(\alpha-i\pi\frac{s-k}{b}\R)\xi}\R)}=2b \frac{\sin\L((s-k)\pi+ib\alpha\R)}{(s-k)\pi+ib\alpha}
\label{eq:26}
\end{equation}
\vspace{0mm}
\noindent By inserting eq.\eqref{eq:25} in eq.\eqref{eq:19}, particularized for $\gamma_k=\rho+i k\Delta\eta$, a system of (complex) algebraic equation is readily found in the unknown variables $X(\gamma_k)$. Once $X(\gamma_k)$ is found for a particular value of $\rho$, the function $x(t)$ is reconstructed by eq.\eqref{eq:15}.

At this stage some comments on the choice of the initial value of $\rho$ are necessary. The value of $\rho$ should belongs to the fundamental strip of both $f(t)$ and $x(t)$. Excluding the case in which $f(t)$ is infinity in zero, the lower limit of its fundamental strip is less than zero. Moreover, if $f(t)\equiv 0$ for $t>t_{\text{max}}$, then the upper limit is certainty $\rho\to\infty$. The limitation $f(t)\equiv 0 \; \forall t>t_{\text{max}}$ is not a problem because of the causal properties of the system at hand, which imply that the response is only influenced by the forcing function in the time interval previous to the instant we are interested in. The lower limit of the fundamental strip of the response $x(t)$, depends on the derivative in zero of $x(t)$ (the system is quiescent up to $t=0$). In particular, it follows from the initial conditions \eqref{eq:17b} that, for $t\to 0$, $x(t)\sim t^n$ so that the lower limit is less than zero and equal to $-n$. In order to know the upper limit of the fundamental strip of $x(t)$, we should know its trend for $t\to\infty$ which, a priori, is not known. We could however guess such a limit to be infinite because this is the value in the case $\lambda=0$ or when the fractional derivative term lacks on the left hand side of eq.\eqref{eq:17a}.

\section{Applications}
In this section, numerical applications of the method outlined in section 3 for the equation \eqref{eq:19}, as well as its generalization to a fractional differential equation with two fractional derivatives of different order, are reported. The method is however general and can be applied for the solution of equations with whatever number of fractional operators. The occurrence of equations with multiple fractional operators, of different order, is peculiar of a wide range of phenomena. Such an example is the Einstein-Smoluchowski equation, ruling the time evolution of a non linear system driven by a general alpha-stable white noise \cite{Cottone2011}. Further examples can be found in acoustical wave propagation in complex media, such as biological tissues, where the resulting fractional wave equations allow to recover the observed power law absorption and dispersion laws \cite{Treeby2010,Prieur2011}, or in diffusion-wave phenomena in general \cite{Mainardi1996}, etc. In order to validate the method, the results are compared with the known solutions in term of the Mittag-Leffler function.

In both cases, the function $x(t)$ has been restored by evaluating the discretized inverse Mellin transform \eqref{eq:15} on the line $\rho=\rho_1=0.5$ of the complex plane, considering a cutoff $\bar{\eta}=200$ and a sampling step $\Delta\eta=0.5$. Evaluating the eq.\eqref{eq:19} in the points $\gamma_k=\rho_1+ik\Delta\eta$ of the complex plane (with $-m\leq k\leq m$; $m=\bar{\eta}/\Delta\eta$), and taking advantage of eq.\eqref{eq:25}, we obtain the system of $2m+1$ equations in the $2m+1$ unknown $X\L(\rho_1+k\Delta\eta\R)$
\vspace{0mm}
\begin{equation}
   \mathbf{M} \mathbf{X}=\mathbf{F}
\label{eq:27}
\end{equation}
\vspace{0mm}
\noindent in which the matrix and vector elements are defined as
\vspace{0mm}
\begin{subequations}
\label{eq:28}
\begin{equation}
   M_{kj}=\frac{1}{2b}\L(C_k^\alpha a_{k-m-1,j-m-1}(\alpha)+\lambda \delta_{kj} (\beta)\R)
\label{eq:28a}
\end{equation}
\begin{equation}
   X_{k}=X(\rho+i(k-m-1)\Delta\eta)
\label{eq:28b}
\end{equation}
\begin{equation}
   F_{k}=F(\rho+i (k-m-1)\Delta\eta)
\label{eq:28c}
\end{equation}
\end{subequations}
\vspace{0mm}
\noindent where $k,j=1\div 2m+1$ and $C_k^{\alpha}=C\L(\rho+i(k-m-1)\Delta\eta,\alpha\R)$. In the figure 1, the solutions obtained for the cases $\lambda=1$ and the forcing function
\vspace{0mm}
\begin{equation}
f(t) =
  \begin{cases}
   \sin{t} & 0\leq t\leq 2\pi \\
   0       & \text{otherwise}
  \end{cases}
\label{eq:29}
\end{equation}
\vspace{0mm}
\noindent is reported, for $\alpha=0.2, 0.5, 0.8, 1.2, 1.5, 1.8$, and compared to the solution in term of Mittag-Leffler function
\vspace{0mm}
\begin{equation}
    x(t)=\int_0^t{\L(t-\tau\R)^{\alpha-1}E_{\alpha,\alpha}\L(-\lambda \L(t-\tau\R)^\alpha\R)f(\tau)\, d\tau}
\label{eq:30}
\end{equation}
\begin{figure}[!h]\centering
\includegraphics[width=10cm]
{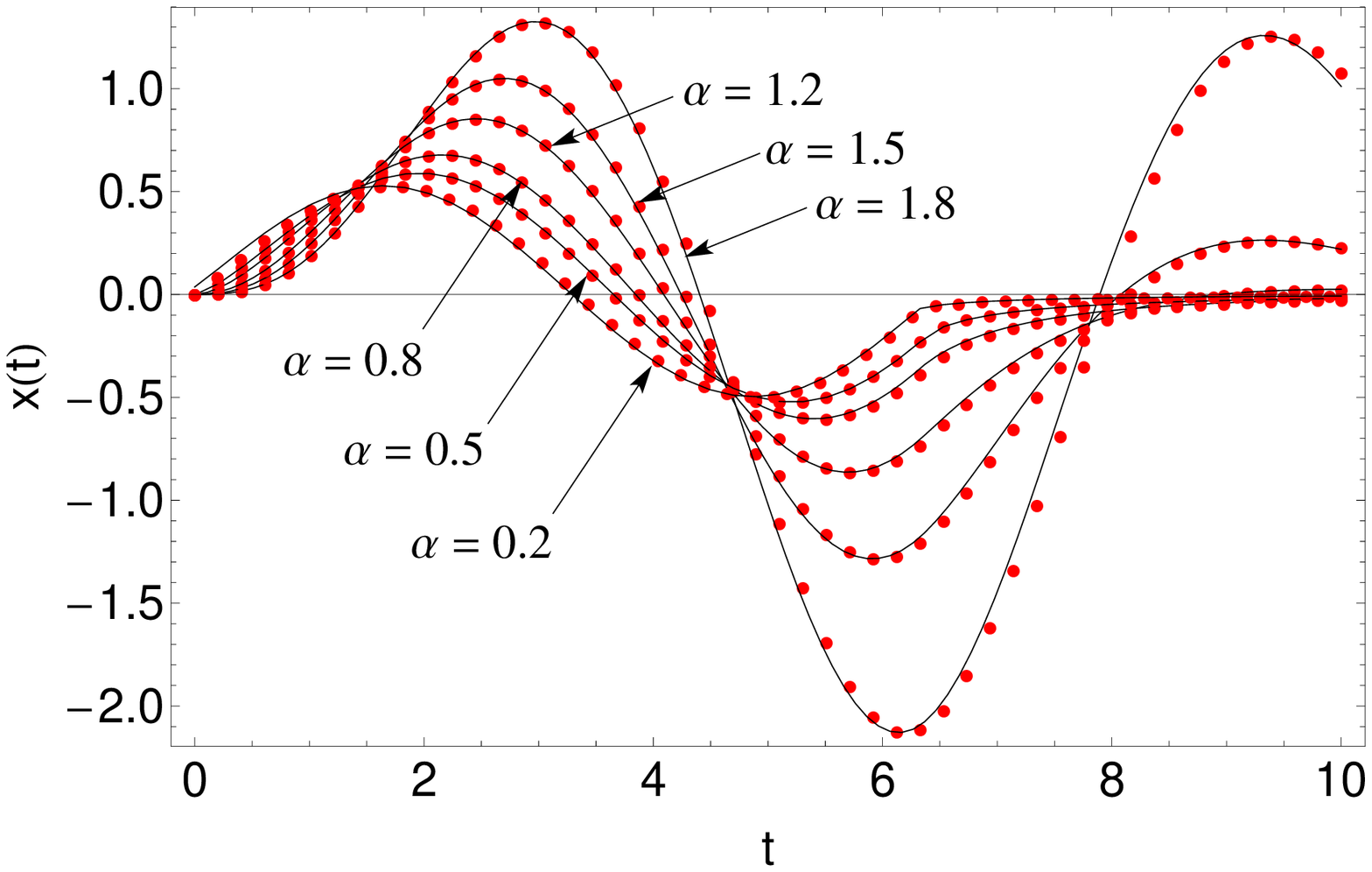}
\caption{Comparison of the solution of eq.\eqref{eq:19} in terms of Mittaf-Leffler function (continuous line) and the solution obtained by \eqref{eq:15}. The value $\lambda=1$ and the forcing function defined in the text have been used.}
\label{fig:1}
\end{figure}
\noindent Let us now consider the fractional differential equation with two fractional derivatives of order $\alpha>0,\; \beta>0$:
\vspace{0mm}
\begin{equation}
   \Dc{x}{\alpha}+\lambda \Dc{x}{\beta}=f(t)
\label{eq:31}
\end{equation}
\vspace{0mm}
\noindent with the same initial conditions eq.\eqref{eq:17b}. Following the same reasoning as above, we obtain the system of $2m+1$ equations
\vspace{0mm}
\begin{equation}
   \mathbf{N} \mathbf{X}=\mathbf{F}
\label{eq:32}
\end{equation}
\vspace{0mm}
\noindent where the matrix elements are now defined as
\vspace{0mm}
\begin{equation}
   N_{kj}=\frac{1}{2b}\L(C_k^\alpha a_{k-m-1,j-m-1}(\alpha)+\lambda C_k^\beta a_{k-m-1,j-m-1}(\beta)\R)
\label{eq:33}
\end{equation}
\vspace{0mm}
\newpage
\noindent with the same meaning of the symbols. In figure 2 is shown the function $x(t)$ restored by \eqref{eq:15}, for the cases $\alpha=0.3$, $\beta=0.5$ and $\alpha=0.2$, $\beta=1.3$, compared to the exact solution
\vspace{0mm}
\begin{equation}
    x(t)=\int_0^t{\L(t-\tau\R)^{\alpha-1}E_{\alpha-\beta,\alpha}\L(-\lambda \L(t-\tau\R)^{\alpha-\beta}\R)f(\tau)\, d\tau}
\label{eq:34}
\end{equation}
\vspace{0mm}
\noindent The value $\lambda=1$ and the forcing function \eqref{eq:29} have been used.
\vspace{0mm}
\begin{figure}[!h]\centering
\includegraphics[width=10cm]
{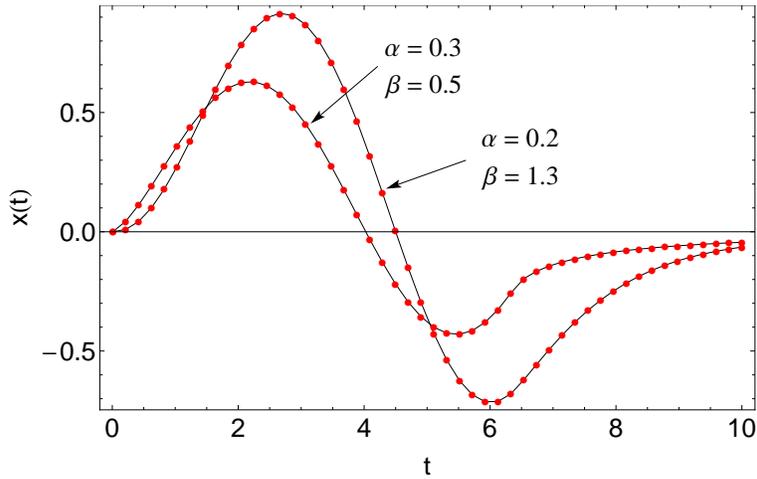}
\caption{Comparison of the solution of eq.\eqref{eq:31} in terms of Mittaf-Leffler function (continuous line) and the solution obtained by \eqref{eq:15}. The value $\lambda=1$ and the forcing function defined in the text have been used.}
\label{fig:2}
\end{figure}
\vspace{0mm}

\vspace{2mm}\noindent For the more general equation
\vspace{0mm}
\begin{equation}
   \sum_k{\lambda_k \Dc{x}{\alpha_k}}=f(t)
\label{eq:35}
\end{equation}
\vspace{0mm}
the solution $x(t)$ may be found in the Mellin domain, in the same fashion as eq.\eqref{eq:32}, with the matrix $\mathbf{N}$ is generalized as
\vspace{0mm}
\begin{equation}
   N_{kj}=\frac{1}{2b}\sum_k{\lambda_k C_k^\alpha a_{k-m-1,j-m-1}(\alpha_j)}
\label{eq:36}
\end{equation}

\newpage
\section{Conclusions}
A method of solution for multi-order fractional differential equations, by using the theory of complex Mellin transform, has been presented. It has been shown that, once the Mellin transform of a function is known in a sufficient number of points of the complex plane $\gamma$, with a fixed value of $\rho=\text{Re}[\gamma]$, the corresponding values of the Mellin transform for any other value $\rho$, say $\rho_j$, may be easily calculated as a linear combination of the former ones (provided $\rho_j$ belongs to the fundamental strip). Taking advantage of this fact, the problem related to the shift of the real part of the argument of the Mellin transform of a fractional operator is overcame, and the fractional differential equation may be easily solved, in the Mellin domain, in terms of a linear set of algebraic equations. The inverse Mellin integral gives, finally, the searched solution in the whole time domain $[0,\infty]$, without invoking Mittag-Leffler functions. The method has been validated by solving two particular fractional differential equations, respectively with one and two fractional operators of different order. The numerical results obtained have been compared with the exact analytic solutions, showing an excellent agreement, regardless of the order of the operators.
The method is robust, computationally efficient, and may be easily implemented for any order and for any number of the fractional derivatives (or integrals) present in the fractional differential equations.

\newpage

\appendix

\section{Fractional calculus}
In this appendix some few preliminaries on fractional calculus will be briefly summarized for introducing the symbols used in the paper. Interested readers are referred to \cite{kilbas2006,podlubny1998,samko1993} for both theory and applications of fractional calculus. Let us start with the definition of Riemann-Liouville (RL) fractional integrals and derivatives. The RL fractional integrals are defined as
\begin{subequations}
	\begin{eqnarray}
			\left( {I_{a_ +  }^\alpha  f} \right)\left( t \right)  &\stackrel{def}{=}  \frac{1}{{\Gamma \left( \alpha  \right)}}\int\limits_a^t {\frac{{f\left( \xi  \right)}}{{\left( {t - \xi } \right)^{1 - \alpha } }}\,d\xi} \;\;\;\;\;\;\;(\alpha \in \mathbb{R},\,t>a) \\
			\left( {I_{b_ -  }^\alpha  f} \right)\left( t \right)   &\stackrel{def}{=} \frac{1}{{\Gamma \left( \alpha  \right)}}\int\limits_t^b {\frac{{f\left( \xi  \right)}}{{\left( {\xi-t } \right)^{1 - \alpha } }}\,d\xi} \;\;\;\;\;\;\; \L(\alpha\in \mathbb{R},\,t<b\R)
	\end{eqnarray}
	\label{eq:A1}
\end{subequations}
where $\Gamma \left( \alpha  \right)$ is the Euler Gamma function. The $\left( {I_{a_ +  }^\alpha  f} \right)\left( x \right)$ and $\left( {I_{b_ -  }^\alpha  f} \right)\left( x \right)$ operators are usually termed as left RL and right RL fractional integrals. Extension to unbounded domain may be simply obtained by letting $a=-\infty$, $b=+\infty$. In the latter case we simply denote them as $\left( {I_{+}^\alpha  f} \right)\left( t \right)$ and $\left( {I_{-}^\alpha  f} \right)\left( t \right)$. The definitions in eq.\eqref{eq:A1} remain meaningful for order $\alpha\in\mathbb{C}$.\\
The RL fractional derivative is defined as
\begin{subequations}
\label{eq:A2}
    \begin{align}
      \left( {D_{a_+}^\alpha  f} \right)\left( t \right) &\stackrel{def}{=}\frac{1}{{\Gamma \left( {n - \alpha } \right)}}\frac{{d^n }}{{dt^n }}\int\limits_a^t {\frac{{f\left( \xi  \right)}}{{\left( {t - \xi } \right)^{\alpha  - n + 1} }}d\xi }\;\; \L(\alpha\in \mathbb{R}, \,t>a\R) \label{eq:A2a} \\
       \left( {D_{b_-}^\alpha  f} \right)\left( t \right) &\stackrel{def}{=}\frac{(-1)^n}{{\Gamma \left( {n - \alpha } \right)}}\frac{{d^n }}{{dt^n }}\int\limits_t^b {\frac{{f\left( \xi  \right)}}{{\left( {t - \xi } \right)^{\alpha  - n + 1} }}d\xi } \;\;\; (\alpha\in \mathbb{R}, \,t<b) \label{eq:A2b}
    \end{align}
  \end{subequations}
where $n=[\alpha]+1 \in \mathbb{N}$ and $[\alpha]$ means integer part of $\alpha$. In the case in which $\alpha=\rho+i\eta\in \mathbb{C}$, eqs.\eqref{eq:A2} hold true posing $n=[\rho]+1$. As well as the fractional integral operators, the RL fractional derivatives could be extended to unbounded domain. We denote $\left( {D_{+}^\alpha  f} \right)\left( x \right)$ and $\left( {D_{-}^\alpha  f} \right)\left( x \right)$ the left and right RL fractional derivatives, respectively.

A definition of fractional derivative useful dealing with physical problems is the Caputo fractional derivative, defined, in its right sided form, as
\begin{equation}
	\Dc{f}{t} \mathop  = \limits^{def} \frac{1}{\Gamma\L(\alpha-n\R)}\int_{a}^{t}{\frac{f^{(n)}(\tau)}{\L(t-\tau\R)^{\alpha+1-n}}\,d\tau} \quad (\alpha\in \mathbb{R}, \,t>a)
\label{eq:A3}
\end{equation}
where again $n=[\alpha]+1 \in \mathbb{N}$. It coalesces to the RL definition for $a\to -\infty$. The main advantage of this definition is that the initial conditions for differential equations of fractional order containing Caputo derivatives, takes the same form as for integer order differential equations, i.e. are expressed in terms of integer order derivatives of the unknown function at the lower limit $t=a$.

Observation of eqs. \eqref{eq:A1}, \eqref{eq:A2} and \eqref{eq:A3} shows that RL fractional operators are convolution integrals with a kernel of power law type. All the rules of ordinary derivatives still apply, including Leibniz's rule and integration by parts. Moreover, Fourier and Laplace transform of RL fractional integrals and derivatives behaves in a quite simple way as in the case of ordinary derivatives and integrals. Such an example, defining the Fourier transform as

\begin{equation}
	{\mathcal F}\left[{f(t)},\omega \right]=\int\limits_{ - \infty }^{+\infty}  {\exp\L({i\omega t}\R)f(t)dt}
\label{eq:A4}
\end{equation}
the following well known relation of the Fourier transform of $d^nf(x)/dt^n$ holds
\begin{equation}
	{\cal F}\left[{\frac{d^n f}{dt^n}},\omega \right]=\int\limits_{ - \infty }^{+\infty}  {\exp\L({i\omega t}\R)\frac{d^n f}{dt^n}dt} =\left(-i \omega\right)^n F(\omega)
\label{eq:A5}
\end{equation}
where $F(\omega)$ is the Fourier transform of $f(t)$. Now, since integration by part remain valid for fractional operators, then it may be easily demonstrated that
\begin{subequations}
\label{eq:A6}
\begin{equation}
   {\cal F} \left\{ {\left( {D_ \pm ^\alpha  f} \right)\left( t \right),\omega } \right\} = \int\limits_{ - \infty }^{+\infty}  {\exp\L({i\omega t}\R) \left( {D_ \pm ^\alpha  f} \right)\left( t \right)dt}  = \left( { \mp i\,\omega } 											\right)^\alpha  F\left( \omega  \right)
\label{eq:A6a}
\end{equation}
\begin{equation}
   {\cal F} \left\{ {\left( {I_ \pm ^\alpha  f} \right)\left( t \right),\omega } \right\} = \int\limits_{ - \infty }^{+\infty}  {\exp\L({i\omega t}\R) \left( {I_ \pm ^\alpha  f} \right)\left( t \right)dt}  = \left( { \mp i\,\omega } 												\right)^{-\alpha}  F\left( \omega  \right)
\label{eq:A6b}
\end{equation}
\end{subequations}
with $\alpha \in \mathbb{C}; \text{Re}(\alpha)=\rho>0$. Analogously, defining the Laplace transform as 
\begin{equation}
	{\mathcal L}\left[{f(t)},s \right]=\int\limits_{0}^{\infty}  {\exp\L({-st}\R)f(t)dt}
\label{eq:A7}
\end{equation}
with $s\in\mathbb{C}$, it can be shown that
\begin{subequations}
	\label{eq:A8}
	\begin{align}
			&{\cal L} \left\{ \L(D_{a^+}^{\alpha}\,f\R)(t);\,s \right\} = s^\alpha F(s)-\sum_{k=0}^{n-1}{s^k\L[D_{a^+}^{\alpha-k-1}f\R]_{t=a}} \;\;\;\;\;\; \L(\alpha\in\mathbb{R}\R) \label{eq:A8a}\\
			&{\cal L} \left\{\L(\prescript{C}{}{\mathbf{D}}_{a^+}^{\alpha} f \R)\L(t\R);\,s \right\} = s^\alpha F(s)-\sum_{k=0}^{n-1}{s^{\alpha-k-1}f^{(k)}(a)} \qquad \L(\alpha\in\mathbb{R}\R) \label{eq:A8b}
	\end{align}
\end{subequations}
where $n=[\alpha]+1 \; \in \; \mathbb{N}$. From eqs.\eqref{eq:A8a} and \eqref{eq:A8b}, the similarity with the integer order differential operators is evident.
\newpage
\bibliography{FractionalOperatorsBiblio}
\bibliographystyle{unsrtnat}
\include{FractionalOperatorsBiblio}

\end{document}